\documentclass[draft,10pt]{article}
\usepackage{amssymb,latexsym,amsmath,array}
\usepackage{makeidx,amsfonts}
\usepackage{longtable,graphicx,ifthen,calc}

\def\qed{\hfill $\square$}
\def\squar{\vbox{\hrule\hbox{\vrule height 6pt \hskip
6pt\vrule}\hrule}}

\def\qed{\hfill $\squar$}
\def\squar{\vbox{\hrule\hbox{\vrule height 6pt \hskip
6pt\vrule}\hrule}}

\begin{document}

\begin{center}
{\LARGE \bf Periodic solutions for planar autonomous systems \\ with nonsmooth periodic perturbations }\\
\bigskip

{\large  Oleg Makarenkov$\,^*$, Paolo Nistri$\,^\dagger$ }\\

\vskip0.2truecm
Dept. of Mathematics, Voronezh State University, Voronezh, Russia \\
e-mail: omakarenkov@kma.vsu.ru\\

\vskip0.2truecm
Dip. di Ingegneria dell' Informazione, Universit\`a di Siena, 53100 Siena, Italy\\
e-mail: pnistri@dii.unisi.it\\
\end{center}

\def\thefootnote{\fnsymbol{footnote}}
\footnotetext[1]{Supported by a President of Russian Federation
Fellowship for Scientific Training Abroad, by the Grant VZ-010 of
RF Ministry of Education and U.S.CRDF, and by RFBR Grants
07-01-00035, 06-01-72552, 05-01-00100. }
\def\thefootnote{\fnsymbol{footnote}}
\footnotetext[2]{Supported by the national research project PRIN
``Control, Optimization and Stability of Nonlinear Systems:
Geometric and Topological Methods''. Corresponding author.}



\begin{abstract}
\noindent In this paper we consider a class of planar autonomous
systems having an isolated limit cycle $x_0$ of smallest period
$T>0$ such that the associated linearized system around it has
only one characteristic multiplier with absolute value 1. We
consider two functions, defined by means of the eigenfunctions of
the adjoint of the linearized system, and we formulate conditions
in terms of them in order to have the existence of two
geometrically distinct families of $T-$periodic solutions of the
autonomous system when it is perturbed by nonsmooth $T-$periodic
nonlinear terms of small amplitude. We also show the convergence
of these periodic solutions to $x_0$ as the perturbation
disappears and we provide an estimation of the rate of
convergence. The employed methods are mainly based on the theory
of topological degree and its properties that allow less
regularity on the data than that required by the approach,
commonly employed in the existing literature on this subject,
based on various versions of the implicit function theorem.

\end{abstract}


\vskip0.5cm\noindent {\bf 1. Introduction}

\noindent Loud in \cite{loud} provided conditions under which the
perturbed system of ordinary differential equations
\begin{equation}\label{ps}
  \dot x=\psi(x)+\varepsilon\phi(t,x,\varepsilon),
\end{equation}
where
\begin{equation}\label{reg}
\psi\in C^2(\mathbb{R}^n,\mathbb{R}^n),\ \ \phi\in
C^1(\mathbb{R}\times\mathbb{R}^n\times[0,1],\mathbb{R}^n)
\end{equation} and $\phi$ is $T$-periodic with respect to time, has, for
sufficiently small $\varepsilon>0,$ a $T$-periodic solution which
tends to a $T$-periodic limit cycle $x_0$ of the unperturbed
system
\begin{equation}\label{np}
  \dot x=\psi(x)
\end{equation}
as $\varepsilon\to 0.$ The limit cycle $x_0$ satisfies the
property that the linearized system
\begin{equation}\label{ls}
  \dot y=\psi'(x_0(t))y
\end{equation}
has only one characteristic multiplier with absolute value 1. Here
and in the following by $C^i(\mathbb{R}^m,\mathbb{R}^n)$ we denote
the vector space of all continuous functions acting from
$\mathbb{R}^m$ to $\mathbb{R}^n$ having $i$-th continuous
derivatives. The main tool employed by Loud is
the following, so-called bifurcation, function
\begin{equation}\label{f0}
  f_0(\theta)=\int_0^T
  \left<z_0(\tau),\phi(\tau-\theta,x_0(\tau),0)\right>d\tau,
\end{equation}
where $z_0$ is a $T$-periodic solution of the adjoint system of
(\ref{ls})
\begin{equation}\label{sopr}
  \dot z=-\left(\psi'(x_0(t))\right)^* z,
\end{equation}
here $A^*$ denotes the transpose of the matrix $A$. Specifically,
Lemma~2 in \cite{loud} states that in order that system (\ref{ps})
has a $T$-periodic solution $x_\varepsilon$ such that
\begin{equation}\label{conv_i}
  x_\varepsilon(t-\theta_0)\to x_0(t){\ \rm as\ }\varepsilon\to 0
\end{equation}
it is necessary that $\theta_0\in\mathbb{R}$ be a zero of the
equation
\begin{equation}\label{eqf0}
f_0(\theta)=0. \end{equation} If (\ref{eqf0}) is satisfied for
some $\theta=\theta_0$ and $f_0'(\theta_0)\not=0,$ i.e. $\theta_0$
is simple, then by (\cite{loud}, Theorem~1) for all sufficiently
small $\varepsilon>0$ system (\ref{ps}) possesses a $T$-periodic
solution $x_\varepsilon$ satisfying
\begin{equation}\label{conv1}
  \|x_\varepsilon(t-\theta_0)- x_0(t)\|\le \varepsilon M,
\end{equation}
where $M>0$ is a constant. These results are also consequences of
general results stated by Malkin in \cite{mal}.

The function $f_0$ has been widely employed to treat different
problems concerning periodic solutions of system (\ref{ps}) with
$\varepsilon>0$ small. We quote in the sequel some papers from the
relevant bibliography devoted to this subject. In \cite{loud} Loud
also  considered the case when (\ref{eqf0}) is identically
satisfied, i.e. $f_0(\theta)=0$ for any $\theta\in[0,T],$ to treat
this case he introduced a new function which plays the role of $f_0$
and he showed that if this function has a simple zero $\theta_0$
then there exists a family of $T$-periodic solutions to (\ref{ps})
satisfying (\ref{conv_i}) (see also \cite{lei}). Moreover in
\cite{loud} it is also considered the case when $\theta_0$ is not a
simple zero of $f_0,$ and the problem of the existence of
$T$-periodic solutions to (\ref{ps}) is associated with the problem
of the existence of roots of a certain quadratic equation. The case
when the limit cycle $x_0$ of system (\ref{np}) is not isolated, in
particular, when the unperturbed system is Hamiltonian, and the case
when the characteristic multiplier of system (\ref{ls}) is not
simple have been considered by many authors. If system (\ref{np}) is
not necessary autonomous and it has a multi-parameterized family of
$T$-periodic solutions, then existence of $T$-periodic solutions of
the perturbed system satisfying (\ref{conv}) was proved by Malkin
\cite{mal}. Melnikov \cite{mel} treated the case when the limit
cycle is not isolated and the limit cycles near $x_0$ are of
different periods (see also Loud \cite{loud1} and Kac \cite{kac})
and he showed that the simple zeros of suitably defined bifurcation
functions $f_{m,n},\ m,n\in\mathbb{N},$ called Melnikov subharmonic
functions, generate periodic solutions in a neighborhood of $x_0$
whose periods are in $m:n$ ratio with respect to the periods of the
perturbation term. Finally, Rhouma and Chicone \cite{chic} have
considered the case when 1 is not a simple multiplier of the
linearized system, to deal with the problem of existence of
$T$-periodic solutions  they introduced a new two variables
bifurcation function $f_0$ whose simple zeros determine families of
$T$-periodic solutions satisfying (\ref{conv_i}).

These theoretical results have been then developed in different
directions: Hausrath and Man\'asevich \cite{man}, (see also
\cite{man1}), found a class of $T$-periodic perturbations $\phi$ for
which the subharmonic Melnikov function $f_{1,1}$ has at least two
simple zeros, obtaining the existence of at least two families of
$T$-periodic solutions to (\ref{ps}) satisfying (\ref{conv_i}).
Makarenkov in \cite{mak} provided useful formulas to calculate
simple zeros of Malkin's bifurcation function in case when the
function $\phi$ is sinusoidal in time. Tkhai \cite{tkh} and Lazer
\cite{laz} developed Malkin's and Melnikov's approaches respectively
to study the existence of periodic solutions to (\ref{ps})
satisfying (\ref{conv_i}) and possessing some additional symmetry
properties that represent relevant features in the applications.
Farkas in \cite{far} investigated the existence of the so-called
$D$-periodic solutions to (\ref{ps}) which are not necessarily
periodic but having periodic derivative. Greenspan and Holmes in
\cite{green} and  Guckenheimer and Holmes in \cite{guck} applied the
method of subharmonic Melnikov's functions to a variety of practical
problems, a number of applications of Malkin's bifurcation function
can be found in the book of Blekhman \cite{ble}.

In all the previous papers, to show the existence of $T$-periodic
solutions for $\varepsilon>0$ small, several formulations of the
implicit function theorem have been employed. Therefore, condition
(\ref{reg}) is the common assumption of these papers (sometimes it
is even required more regularity on $\psi$ and $\phi$). The
persistence of the limit cycle $x_0$ under less restrictive
regularity assumptions than (\ref{reg}) is studied only for the
cases when system (\ref{np}) is linear, in this case the modified
averaging methods developed by Mitropol'sksii \cite{mit} and
Samoylenko \cite{sam} can be applied as well as the coincidence
degree theory introduced by Mawhin, see, for instance,
(\cite{mawhin}, Theorem~IV.13); Hamiltonian, see M.~Henrard and
F.~Zanolin \cite{zan}; or piecewise differentiable, see
Kolovski\u\i\ \cite{kol} and \v Ste\u\i nberg \cite{ste}.

In the present paper we assume that the linearized system (\ref{ls})
has only one characteristic multiplier equal to $1$ and
\begin{equation}\label{reg1}
\psi\in C^1(\mathbb{R}^2,\mathbb{R}^2),\ \ \phi\in
C(\mathbb{R}\times\mathbb{R}^2\times[0,1],\mathbb{R}^2).
\end{equation}
By combining the function $f_0$ with the analogously defined
function
$$
  f_1(\theta,s)=\int_{s-T}^s \left<
  z_1(\tau),\phi(\tau-\theta,x_0(\tau),0)\right>d\tau,
$$
where $z_1$ is an eigenfunction of system (\ref{sopr})
corresponding to the characteristic multiplier $\rho_*\not=1,$ we
give conditions in Theorem~3 for the existence of $T$-periodic
solutions to (\ref{ps}) satisfying (\ref{conv_i}). Although, as we
have mentioned before, in many papers it was proved the existence
of two or more families of $T$-periodic solutions to (\ref{ps})
converging to $x_0$ in the sense of (\ref{conv_i}), it was not
guaranteed that these families do not coincide geometrically,
namely if one is just a shift in time of the other. In this paper
our results ensure the existence of at least two geometrically
distinct families of $T$-periodic solutions to (\ref{ps})
satisfying (\ref{conv_i}). Moreover, since property (\ref{conv1})
is a consequence of the application of the implicit function
theorem it is not anymore guaranteed under our conditions
(\ref{reg1}). However, we will show in Theorem~1 that under
conditions (\ref{reg1}) the following property holds
\begin{equation}\label{conv2i}
  \varepsilon M_1 |f_1(\theta_0,t)| \le \left\|x_\varepsilon(t-\theta_\varepsilon(t))-x_0(t)
  \right\|+o(\varepsilon)\le \varepsilon M_2
  |f_1(\theta_0,t)|\quad {\rm for\ any\
  }\varepsilon\in(0,\varepsilon_0)\  {\rm and\ any\
  } t\in[0,T],
\end{equation}
where $0<M_1<M_2$ and $\theta_\varepsilon(t)\to \theta_0$ as
$\varepsilon\to 0$ uniformly with respect to $t\in[0,T].$ The
introduction of the function $f_1,$ as shown by (\ref{conv2i}) and
Corollaries 1 and 2 of this paper, gives a new qualitative
information about the convergence (\ref{conv_i}) with respect to
(\ref{conv1}) and it is a contribution to the problem posed by
Hale and T\'boas in \cite{hale} concerning the behavior of the
periodic solutions of a second order periodically perturbed
autonomous system when the perturbation disappears. We would like
also to remark, that Loud in \cite{loud} provided a precise
information about the way of convergence of $x_\varepsilon$ to
$x_0$ by means of  the representation
$x_\varepsilon(t)=x_0(t+\theta_0)+\varepsilon
y(t+\theta_0)+o(\varepsilon),$ where the function $y$ is a
suitably chosen solution of system (\ref{aux}) of this paper with
$\xi=x_0(0),$ see (\cite{loud}, formulas 1.3 and 2.11).

In order to prove the existence of $T$-periodic solutions to
(\ref{ps}) satisfying (\ref{conv_i}) under assumptions (\ref{reg1})
we make use of the topological degree theory. Specifically, for
$\varepsilon>0,$ we consider  the integral operator
$G_\varepsilon:C([0,T],\mathbb{R}^2)\to C([0,T],\mathbb{R}^2)$ given
by $(G_\varepsilon x)(t)=x(T)+\int_0^t
  \psi(x(\tau))d\tau+\varepsilon
  \int_0^t\phi(\tau,x(\tau),\varepsilon)d\tau,\, t\in [0,T],$
here $C([0,T],\mathbb{R}^2)$ is the Banach space of all the
continuous functions defined on $[0,T]$ with values in
$\mathbb{R}^2$ equipped with the sup-norm. We also consider the
Leray-Schauder degree $d(I-G_\varepsilon,W_U),$ see Brown
(\cite{bro}, \S 9), of the compact vector field $I-G_\varepsilon$
with respect to the open set $W_U=\left\{x\in
C([0,T],\mathbb{R}^2):x(t)\in U,\right.$ for any $\left.
t\in[0,T]\right\},$ where $U$ is an open set of $\mathbb{R}^2.$ We
will provide conditions in terms of the functions $f_0$ and $f_1$
ensuring that $d(I-G_\varepsilon,W_{U_0})\not=
d(I-G_\varepsilon,W_{U_\varepsilon})$ for all $\varepsilon>0$
sufficiently small, where $U_0$ is the interior of the limit cycle
$x_0$ and $U_\varepsilon\subset U_0$ (or $U_0\subset
U_\varepsilon$) is a suitably defined family of sets such that
$U_\varepsilon\to U_0$ as $\varepsilon\to 0.$ To do this we use a
result by Capietto, Mawhin and Zanolin (\cite{maw}, Corollary~1)
which, under our assumptions, states that
$d(I-G_0,W_{U_\varepsilon})=1$ for $\varepsilon>0$ sufficiently
small. Then in Theorem~2, by means of a result due to Kamenskii,
Makarenkov and Nistri (\cite{non}, Theorem~2),  we conclude that
there exists a continuous vector field $F: \mathbb{R}^2\to
\mathbb{R}^2$ with
\begin{equation}
\label{repi}F(x_0(\theta))=f_0(\theta)\dot
x_0(\theta)+f_1(\theta,\theta)y_1(\theta),\quad {\rm for\ any\
  } \theta\in [0,T],
\end{equation}
here $y_1$ denotes the eigenfunction of (\ref{ls}), corresponding to
the characteristic multiplier $\rho\not=1,$ such that for all
$\varepsilon>0$ sufficiently small we have that
$d(I-G_\varepsilon,W_{U_0})=d_B(F,U_0),$ where $d_B(F,U_0)$ is the
Brouwer degree of $F$ on $U_0,$ see e.g. Brown (\cite{bro}, \S 8).
In our case the integer $d_B(F,U_0)$ can be easily calculated since
it is equal to the Poincar\'e index of $x_0$ with respect to the
vector field $F$ multiplied by $+1$ or $-1$ according with the
orientation of $x_0,$ see Lefschetz (\cite{lef}, Ch. IX, \S 4).
Furthermore, as it was observed by Bobylev and Krasnoselskii in
\cite{bob}, for any small neighborhood $B_\delta(\partial W_{U_0})$
of the boundary $\partial W_{U_0}$ of $W_{U_0},$ we have that
$d(I-G_0,B_\delta(\partial W_{U_0}))=0$ and so one cannot directly
apply Leray-Schauder fixed point theorem for studying the existence
of $T$-periodic solutions to (\ref{ps}) satisfying (\ref{conv_i}).

The paper is organized as follows. Theorem~1 of Section~2 states
property (\ref{conv2i}) for the $T$-periodic solutions of system
(\ref{ps}), in Theorem~2 we prove the coincidence degree formula
$d(I-G_\varepsilon,W_{U_0})=d_B(F,U_0)$ for $\varepsilon>0$ small.
Finally, we give the main result of the paper: Theorem~3 which
states the existence of at least two geometrically distinct
families of $T$-periodic solutions to (\ref{ps}) satisfying
(\ref{conv_i}). In Section 3 we provide an example which shows how
formula (\ref{repi}) can be used for the practical calculation of
$d_B(F,U_0).$ In fact, under quite general conditions on $f_0$ and
$f_1,$ we show that if $\phi(t,\xi)=-\phi(t+T/2,\xi)$ for any
$t\in [0,T]$ and $\xi\in \mathbb{R}^2$, then
$d_B(F,U_0)\in\{0,2\}.$ Finally, we outline some methods for
calculating the eigenfunctions $y_1,$ $z_0$ and $z_1.$

\vskip0.5cm \noindent
 {\bf 2. Main results. }

\noindent Through the paper we assume the following condition:

\vskip0.3cm

\noindent ($A_0$)$-$ system (\ref{np}) has a limit cycle $x_0$ with
smallest period $T>0$ and the linearized system (\ref{ls}) has only
one characteristic multiplier equal to $1.$

\vskip0.2cm

\noindent In what follows we provide the notations that we will
use in the proofs of the results of this Section.

\vskip0.2cm\noindent By $y_1$ we denote the eigenfunction of
(\ref{ls}) corresponding to the characteristic multiplier
$\rho\not=1$ (clearly $\dot x_0$ is the eigenfunction of
(\ref{ls}) corresponding to the characteristic multiplier $1$).
Moreover, $z_0$ and $z_1$ will denote the eigenfunctions of
(\ref{sopr}) corresponding to the characteristic multipliers $1$
and $\rho_*\not=1$ respectively. $(a_1,a_2)$ is the matrix whose
columns are the vectors $a_1,a_2\in\mathbb{R}^2,$ $(a_1,a_2)^*$
denotes the transpose of $(a_1,a_2)$, $\Omega(\cdot,t_0,\xi)$ is
the solution of system (\ref{np}) satisfying $x(t_0)=\xi,$
$\Omega'_{\xi}(\cdot,t_0,\xi)$ is the derivative of
$\Omega(\cdot,t_0,\xi)$ with respect to the third variable, $U_0$
is the interior of the limit cycle $x_0$ of system (\ref{np}),
$\partial U_0$ is the boundary of $U_0,$ $[v]_i,$ $i=1,2,$ is the
$i$-th component of vector $v\in\mathbb{R}^2.$ $a\parallel b$ will
indicate that the vectors $a,b\in\mathbb{R}^2$ are parallel,
$a^\bot$ denotes the vector $a\in\mathbb{R}^2$ rotated of $\pi/2$
clockwise and
$\angle(a,b)=\arccos\dfrac{\left|\left<a,b\right>\right|}{\|a\|\cdot\|b\|}$
is the angle between the two vectors $a,b\in\mathbb{R}^2$. By
$o(\varepsilon),$ $\varepsilon>0,$ we will denote a function,
which may depend also on other variables having the property that
$\dfrac{o(\varepsilon)}{\varepsilon}\to 0$ as $\varepsilon\to 0$
uniformly with respect to the other variables when they belong to
a bounded set.

\noindent Finally, let $t, r\in \mathbb{R}$ and let $h(t,r)$ be the
vector of $\mathbb{R}^2$ given by
\begin{equation}\label{h}
 h(t,r)=x_0(t)+r\frac{z_0(t)^\bot}{\|z_0(t)^\bot\|}.
\end{equation}
Define the function $(t,r)\to I(t,r)$ as follows
\begin{equation}\label{I}
  I(t,r)=\Omega(T,0,h(t,r)).
\end{equation}
It is easily seen that, for any $t\in [0,T],$ the curve $r\to
I(t,r)$ intersects the limit cycle $x_0$ at the point
$I(t,0)=x_0(t).$

\vskip0.4truecm \noindent The following theorem states a property
similar to (\ref{conv1}) in the case when the autonomous system
(\ref{np}) is perturbed by nonsmooth functions $\phi.$

\noindent {\bf Theorem 1.} {\it Assume  conditions (\ref{reg1}).
Assume that, for all sufficiently small $\varepsilon>0,$ system
(\ref{ps}) has a $T$-periodic solution $x_\varepsilon$ satisfying
\begin{equation}\label{conv}
x_\varepsilon(t-\theta_0)\to x_0(t){\ \ as\ \ }\varepsilon\to 0,
\end{equation}
for any $t\in [0,T]$, where $\theta_0\in[0,T].$ Then there exist
constants $ 0<M_1<M_2,$ $\varepsilon_0>0$ and $r_0\in (0,1]$ such
that
\begin{equation}\label{conv2}
  \varepsilon M_1 |f_1(\theta_0,t)| \le \left\|x_\varepsilon(t-\theta_\varepsilon(t))-
  x_0(t)\right\|+o(\varepsilon)\le \varepsilon M_2
  |f_1(\theta_0,t)|\quad { for\ any\
  }\varepsilon\in(0,\varepsilon_0)\ { and\ any\ }t\in[0,T],
\end{equation}
where $\theta_\varepsilon(t)\to \theta_0$ as $\varepsilon\to 0$
uniformly with respect to $t\in[0,T],$ and $\
x_\varepsilon(t-\theta_\varepsilon(t))\in I(t,[-r_0,r_0]),$
$t\in[0,T].$}

\vskip0.2cm\noindent To prove Theorem 1 we need some preliminary
lemmas.

\noindent {\bf Lemma 1.} {\it For any $t\in\mathbb{R}$ we have
\begin{equation}\label{sat}
  (\dot x_0(t)\ y_1(t))^*(z_0(t)\ z_1(t))=\left(\begin{array}{cc} \left<\dot x_0(0),z_0(0)\right> & 0\\
  0 & \left<y_1(0),z_1(0)\right> \end{array}\right).
\end{equation} }

\noindent {\bf Proof.}  By Perron's lemma \cite{perron} (see also
Demidovich (\cite{dem}, Sec. III, \S 12) for any $t\in\mathbb{R}$ we
have
$$
(\dot x_0(t)\ y_1(t))^*(z_0(t)\ z_1(t)):=
\left(\begin{array}{ll}\left<\dot x_0(t),z_0(t)\right>
  & \left<\dot x_0(t),z_1(t)\right> \\ \left< y_1(t),z_0(t)\right> &
  \left<y_1(t),z_1(t)\right>\end{array}\right)=\left(\begin{array}{ll}\left<\dot x_0(0),z_0(0)\right>
  & \left<\dot x_0(0),z_1(0)\right> \\ \left< y_1(0),z_0(0)\right> &
  \left<y_1(0),z_1(0)\right>\end{array}\right).
$$
Thus, in particular, $\left<\dot x_0(0),z_1(0)\right>=\left<\dot
x_0(T),z_1(T)\right>.$ On the other hand $\dot x_0(0)=\dot x_0(T)$
and $z_1(T)=\rho_* z_1(0),$ $\rho_*\not=1,$ thus $\left<\dot
x_0(0),z_1(0)\right>=0.$ Analogously, since $y_1(T)=\rho y_1(0),$
$\rho\not=1,$ and $z_0(0)=z_0(T),$ we have that
$\left<y_1(0),z_0(0)\right>=0.$

                                                              \qed

\noindent {\bf Lemma 2.} {\it Under the assumptions of Theorem 1
there exist $r_0\in (0,1]$ and $\alpha_0\in [0,\pi/2)$ such that

\begin{equation}\label{st}
   \angle(I(t,r)-x_0(t),\dot
   x_0(t)^\bot)<\alpha_0 \quad {\rm for\ any\
   }t\in[0,T]\quad {\rm and\ any\
   }r\in[-r_0,r_0].
\end{equation}}

\noindent {\bf Proof.} Assume the contrary, hence there exist
sequences $\{t_n\}_{n\in\mathbb{N}}\subset[0,T],$ $t_n\to t_0$ as
$n\to\infty,$ $\{r_n\}_{n\in\mathbb{N}}\subset(0,1],$ $r_n\to 0$ as
$n\to\infty$ such that
\begin{equation}\label{cw}
   \angle(I(t_n,r_n)-x_0(t_n),\dot
   x_0(t_n)^\bot)\to \pi/2\quad{\rm as\ } n\to\infty.
\end{equation}
We have
\begin{equation}\label{pp}
 I(t_n,r_n)-x_0(t_n)=\Omega\left(T,0,x_0(t_n)+r_n
 \frac{z_0(t_n)^\bot}{\|z_0(t_n)^\bot\|}\right)-x_0(t_n)=
 r_n\Omega'_{\xi}(T,0,x_0(t_n))\frac{z_0(t_n)^\bot}{\|z_0(t_n)^\bot\|}+o(r_n).
\end{equation}
By Theorem 2.1 of \cite{kra} it follows that
$\Omega'_\xi(T,0,h(t,0))=Y(T,t)$ where $Y(\cdot,t)$ is the
fundamental matrix for the system
\begin{equation}\label{lst}
  \dot y(\tau)=\psi'(x_0(\tau+t))y(\tau),
\end{equation}
satisfying $Y(0,t)=I,$ thus $\Omega'_\xi(T,0,h(t,0))y_1(t)=\rho
y_1(t).$ On the other hand from Lemma~1 we have
\begin{equation}\label{fl1}
y_1(t)\parallel z_0(t)^\bot,
\end{equation} therefore
$\Omega'_\xi(T,0,h(t,0))z_0(t)^\bot=\rho z_0(t)^\bot$ and (\ref{pp})
can be rewritten as follows
\begin{equation}\label{pp1}
 I(t_n,r_n)-x_0(t_n)=\rho
 r_n\frac{z_0(t_n)^\bot}{\|z_0(t_n)^\bot\|}+o(r_n).
\end{equation}
Hence
$$
  \angle(I(t_n,r_n)-x_0(t_n),\dot
   x_0(t_n)^\bot)=\arccos\dfrac{\left|\left<\rho
 \dfrac{z_0(t_n)^\bot}{\|z_0(t_n)^\bot\|}+\dfrac{o(r_n)}{r_n},
 \dot x_0(t_n)^\bot\right>\right|}{\left\|\rho
 \dfrac{z_0(t_n)^\bot}{\|z_0(t_n)^\bot\|}+\dfrac{o(r_n)}{r_n}
 \right\|\cdot\left\|\dot x_0(t_n)^\bot\right\|}.
$$
Without loss of generality we may assume
\begin{equation}\label{sign1}
\left<\dot x_0(0),z_0(0)\right>=1
\end{equation}
and so
$$
\angle(I(t_n,r_n)-x_0(t_n),\dot
   x_0(t_n)^\bot)\to\arccos\frac{1}{\|z_0(t_0)^\bot\|\cdot\left\|\dot
   x_0(t_0)^\bot\right\|}\qquad {\rm as\ }n\to\infty
$$
contradicting (\ref{cw}). Therefore there exist $r_0\in(0,1]$ and
$\alpha_0\in[0,\pi/2)$ satisfying (\ref{st}).

                                                               \qed

\noindent We can now prove the following.

\noindent {\bf Lemma 3.} {\it Under the assumptions of Theorem 1
there exists $\varepsilon_0>0$ such that for any $\varepsilon\in
(0,\varepsilon_0)$ and any $t\in [0,T]$ we have
$$
x_\varepsilon(t-\theta_\varepsilon(t))\in I(t,[-r_0,r_0])
$$
where $\theta_\varepsilon(t)=\theta_0-\Delta_\varepsilon(t),$
$\Delta_\varepsilon(t)\in [t-\frac{T}{2}, t+\frac{T}{2}]$ and
$\Delta_\varepsilon(t)\to 0$ as $\varepsilon \to 0$ uniformly in
$t\in [0,T]$. Moreover, there exists $M>0$ such that
\begin{equation}\label{ST}
 \left\|x_\varepsilon(t-\theta_\varepsilon(t))-x_0(t)\right\|\le\varepsilon M\quad\mbox{for\ any\ }
 \varepsilon\in(0,\varepsilon_0) \quad\mbox{and\ any\ } t\in [0,T].
\end{equation}
}

\noindent {\bf Proof.} First of all observe that $r_0>0$ given by
Lemma~2 can be chosen to satisfy
\begin{equation}\label{sta1}
  I(t,[-r_0,r_0])\cap x_0([0,T])=\{x_0(t)\},{\ \rm for\
  any\ }t\in[0,T].
\end{equation}
From (\ref{st}) of Lemma 2 and (\ref{conv}) we have that there
exists $\varepsilon_0>0$ such that $I(t,[-r_0,r_0])\cap
x_\varepsilon([0,T])\not=\emptyset$ for any
$\varepsilon\in(0,\varepsilon_0)$ and any $t\in[0,T].$ Hence, for
any $\varepsilon\in (0,\varepsilon_0)$ and $t\in[0,T]$ there exists
$\Delta_\varepsilon(t)\in\left[t-\frac{T}{2},t+\frac{T}{2}\right]$
such that
\begin{equation}\label{111}
x_\varepsilon(t-\theta_0+\Delta_\varepsilon(t))\in I(t,[-r_0,r_0]).
\end{equation}
We claim that
\begin{equation}\label{convD}
  \Delta_\varepsilon(t)\to 0\ {\rm  as\ } \varepsilon\to 0
\end{equation}
uniformly with respect to $t\in[0,T].$ In fact, assume the contrary,
thus there exist sequences
$\{\varepsilon_n\}_{n\in\mathbb{N}}\subset(0,\varepsilon_0),$
$\varepsilon_n\to 0$ as $n\to\infty,$ and
$\{t_n\}_{n\in\mathbb{N}},$ $t_n\to t_0\in[0,T]$ as $n\to\infty,$
such that $\Delta_{\varepsilon_n}(t_n)\to\Delta_0\not=0$ and
$x_{\varepsilon_n}(t_n-\theta_0+\Delta_{\varepsilon_n}(t_n))\in
I(t_n,[-r_0,r_0]).$ Since from (\ref{sta1}) we have
$I(t_n,[-r_0,r_0])\cap x_0([0,T])=\{x_0(t_n)\}$ then
\begin{equation}\label{w1}
x_{\varepsilon_n}(t_n-\theta_0+\Delta_{\varepsilon_n}(t_n))\to
x_0(t_0)\quad {\rm as\ }n\to\infty.
\end{equation}
Applying (\ref{conv}) we have
\begin{equation}\label{w2}
x_{\varepsilon_n}(t_n-\theta_0+\Delta_{\varepsilon_n}(t_n))\to
x_0(\Delta_0+t_0).
\end{equation}
From (\ref{w1}) and (\ref{w2}) we conclude that
\begin{equation}\label{star1}
x_0(t_0)=x_0(\Delta_0+t_0),
\end{equation}
where $\Delta_0\in\left[t_0-\frac{T}{2},t_0+\frac{T}{2}\right],$
since $T$ is the smallest period of $x_0$ it follows from
(\ref{star1}) that $\Delta_0=0,$ which is a contradiction.

\vskip2mm \noindent Pick any $\tau\in [0,T]$, in what follows we
show that the shifts $t \to\Delta_\varepsilon(t)$ have the
property that the convergence of
$x_\varepsilon(\tau+t-\theta_\varepsilon(t))$ to $x_0(\tau+t)$ is
of order $\varepsilon>0,$ where
$\theta_\varepsilon(t)=\theta_0-\Delta_\varepsilon(t),$ and thus
the claim of Lemma~3 is proved.

\noindent For this consider the change of variables
$\nu_\varepsilon(\tau,t)=\Omega(0,\tau,x_\varepsilon(\tau+t-\theta_\varepsilon(t)))$
in system (\ref{ps}). It is clear that
$x_\varepsilon(\tau+t-\theta_\varepsilon(t))=\Omega(\tau,0,\nu_\varepsilon(\tau,t))$
and so
\begin{equation}\label{ob1}
 \dot x_\varepsilon(\tau+t-\theta_\varepsilon(t))=
 \psi(\Omega(\tau,0,\nu_\varepsilon(\tau,t))+\Omega'_{\xi}(\tau,0,
 \nu_\varepsilon(\tau,t)) (\nu_\varepsilon)'_\tau(\tau,t).
\end{equation}
On the other hand from (\ref{ps}) we have
\begin{equation}\label{ob2}
  \dot
  x_\varepsilon(\tau+t-\theta_\varepsilon(t))=
  \psi(\Omega(\tau,0,\nu_\varepsilon(\tau,t)))+\varepsilon
  \phi(\tau+t-\theta_\varepsilon(t),
  \Omega(\tau,0,\nu_\varepsilon(\tau,t)),\varepsilon).
\end{equation}
From (\ref{ob1}) and (\ref{ob2}) it follows
$$
  (\nu_\varepsilon)'_\tau(\tau,t)=\varepsilon\left(\Omega'_{\xi}
  (\tau,0,\nu_\varepsilon(\tau,t))\right)^{-1}\phi(\tau+t-\theta_\varepsilon(t),\Omega(\tau,0,\nu_\varepsilon(\tau,t)),\varepsilon)
$$
and since
\begin{equation}\label{ob2bis}
  \nu_\varepsilon(0,t)=x_\varepsilon(t-\theta_\varepsilon(t))=x_\varepsilon(T+t-\theta_\varepsilon(t))=
  \Omega(T,0,\nu_\varepsilon(T,t))
\end{equation}
we finally obtain
\begin{equation}\label{ob3}
  \nu_\varepsilon(\tau,t)=\Omega(T,0,\nu_\varepsilon(T,t))+\varepsilon\int_0^\tau
  \left(\Omega'_{\xi}(s,0,\nu_\varepsilon(s,t))\right)^{-1}\phi(s+t -\theta_\varepsilon(t),\Omega(s,0,\nu_\varepsilon(s,t)),\varepsilon)
  ds.
\end{equation}
Since $\nu_\varepsilon(\tau,t)\to \Omega(0,\tau,x_0(\tau+t))=x_0(t)$
as $\varepsilon\to 0$ we can write $\nu_\varepsilon(\tau,t)$ in the
following form
\begin{equation}\label{rep}
  \nu_\varepsilon(\tau,t)=x_0(t)+\varepsilon
  \mu_\varepsilon(\tau,t).
\end{equation}
Subtract $x_0(t)$ from both sides of (\ref{ob3}) obtaining
\begin{eqnarray}
\varepsilon\mu_\varepsilon(\tau,t)&=&\varepsilon\,\Omega'_{\xi}(T,0,x_0(t))
\mu_\varepsilon(T,t)+o(\varepsilon\mu_\varepsilon(T,t))+\nonumber\\
& &+\>
\varepsilon\int_0^\tau\left(\Omega'_{\xi}(s,0,\nu_\varepsilon(s,t))
\right)^{-1}\phi(s+t-\theta_\varepsilon(t),
\Omega(s,0,\nu_\varepsilon(s,t)),\varepsilon)
  ds.\label{ob4}
\end{eqnarray}
Since $x_\varepsilon(t-\theta_\varepsilon(t))\in I(t,[-r_0,r_0])$
then from (\ref{I}) there exists $r_\varepsilon(t)\in[-r_0,r_0]$
such that $x_\varepsilon(t-\theta_\varepsilon(t))=
\Omega(T,0,h(t,r_\varepsilon(t)))$ and by (\ref{h}) we get
\begin{eqnarray*}\varepsilon\mu_\varepsilon(T,t) & =
& \nu_\varepsilon(T,t)-x_0(t)=\Omega(0,T,x_\varepsilon
(t-\theta_\varepsilon(t)))-x_0(t)=\\
& =& \Omega(0,T,\Omega(T,0,h(t,r_\varepsilon(t))))-x_0(t)=
h(t,r_\varepsilon(t))-x_0(t)=r_\varepsilon(t)
\frac{z_0(t)^\bot}{\|z_0(t)^\bot\|}.
\end{eqnarray*}
Therefore $\mu_\varepsilon(T,t)\parallel z_0(t)^\bot$ and by
(\ref{fl1}) we can rewrite (\ref{ob4}) as follows
\begin{equation}
\varepsilon\mu_\varepsilon(\tau,t)=\varepsilon\rho\mu_\varepsilon(T,t)+
o(\varepsilon\mu_\varepsilon(T,t))+
\varepsilon\int_0^\tau\left(\Omega'_{\xi}(s,0,\nu_\varepsilon(s,t))
\right)^{-1}\phi(s+t-\theta_\varepsilon(t),
\Omega(s,0,\nu_\varepsilon(s,t)),\varepsilon)
  ds.\label{bisp}
\end{equation}
We now prove that the functions $(\tau,t)\to
\mu_\varepsilon(\tau,t)$ are uniformly bounded with respect to
$\varepsilon\in(0,\varepsilon_0).$ For this we argue by
contradiction, therefore there exist sequences
$\{\varepsilon_n\}_{n\in\mathbb{N}}\subset(0,1),$ $\varepsilon_n\to
0$ as $n\to\infty,$ $\{\tau_n\}_{n\in\mathbb{N}}\subset[0,T],$
$\tau_n\to \tau_0$ as $n\to\infty$ and
$\{t_n\}_{n\in\mathbb{N}}\subset[0,T],$ $t_n\to t_0$ as
$n\to\infty,$ such that
$\|\mu_{\varepsilon_n}(\tau_n,t_n)\|\to\infty$ as $n\to\infty,$ so
$\|\mu_{\varepsilon_n}(\cdot,t_n)\|_C\to\infty$ as $n\to\infty,$
where $\|\cdot\|_C$ is the usual sup-norm of
$C([0,T],\mathbb{R}^2).$ Let
$q_n(\tau)=\dfrac{\mu_{\varepsilon_n}(\tau,t_n)}{\|\mu_{\varepsilon_n}(\cdot,t_n)\|_C},$
then from (\ref{bisp}) we have
\begin{eqnarray}\label{ff}
 q_n(\tau)&=&\rho
 q_n(T)+\frac{o(\varepsilon_n
 \mu_{\varepsilon_n}(T,t_n))}{\varepsilon_n\|\mu_{\varepsilon_n}
 (\cdot,t_n)\|_C}+\nonumber\\
 &&+\>\frac{1}{\|\mu_{\varepsilon_n}(\cdot,t_n)\|_C}\int_0^\tau
 \left(\Omega'_{\xi}(s,0,\nu_{\varepsilon_n}(s,t_n))
 \right)^{-1}\phi(s+t_n-\theta_{\varepsilon_n}(t_n),\Omega(s,0,\nu_{\varepsilon_n}(s,t_n)),\varepsilon_n)
  ds.
\end{eqnarray}
By definition the set of continuous functions $A=\{q_n,\
n\in\mathbb{N}\},$ is bounded and, as it is easy to see from
(\ref{ff}), $A$ is also equicontinuous. Therefore, by the
Ascoli-Arzela Theorem, see e.g. (\cite{bro}, Theorem~2.3), we may
assume without loss of generality that the sequence
$\{q_n\}_{n\in\mathbb{N}}$ is converging. Let
$q_0=\lim_{n\to\infty} q_n,$ from (\ref{ff}) we may conclude that
\begin{equation}\label{ob5}
 q_0(\tau)=\rho q_0(T).
\end{equation}
By (\ref{ob5}) it follows that $q_0$ is a constant function, thus
being $\rho\not=1$  we have $q_0=0.$ On the other hand, by the
definition of $q_n,$ we have that $\|q_0\|_C=1.$ This contradiction
shows the uniform boundedness of the functions $\mu_\varepsilon$
with respect to $\varepsilon\in (0, \varepsilon_0).$ On the other
hand from (\ref{ob2bis}) and (\ref{rep})we have that
\begin{equation}\label{wehav}
x_\varepsilon(t-\theta_\varepsilon(t))-x_0(t)=\varepsilon\mu_\varepsilon(0,t),
\end{equation}
and thus the proof is complete.

                                                          \qed

\noindent {\bf Proof of Theorem 1.} We have to prove (\ref{conv2})
with $t\to\theta_\varepsilon(t)$ as given in Lemma~3. For this, by
Lemma~1, we can represent
$x_\varepsilon(\tau+t-\theta_\varepsilon(t))-x_0(\tau+t)$ as
follows
\begin{equation}\label{fw}
 x_\varepsilon(\tau+t-\theta_\varepsilon(t))-x_0(\tau+t)=\varepsilon
  a_\varepsilon(\tau,t)\dot x_0(\tau+t)+\varepsilon
  b_\varepsilon(\tau,t)y_1(\tau+t),
\end{equation}
where
\begin{eqnarray}
& &\varepsilon
a_\varepsilon(\tau,t)=\left<z_0(\tau+t),x_\varepsilon(\tau+t-\theta_\varepsilon(t))-x_0(\tau+t)\right>
\qquad\mbox{and}\nonumber\\
& &\varepsilon
b_\varepsilon(\tau,t)=\left<z_1(\tau+t),x_\varepsilon(\tau+t-\theta_\varepsilon(t))-x_0(\tau+t)\right>.\label{b}
\end{eqnarray}
By Lemma 1 we have that $\left<\dot x_0(t),z_1(t)\right>=0,$ for any
$t\in[0,T]$, and so $\dot x_0(t)^\bot=k\dfrac{\|\dot
x_0(t)^\bot\|}{\|z_1(t)\|}z_1(t),$ where $k=+1$ or $k=-1.$ Therefore
\begin{equation}\label{form}
\left<\dot
x_0(\tau+t)^\bot,x_\varepsilon(\tau+t-\theta_\varepsilon(t))-x_0(\tau+t)\right>=\varepsilon
b_\varepsilon(\tau,t)k\frac{\|\dot
x_0(\tau+t)^\bot\|}{\|z_1(\tau+t)\|}.
\end{equation}
We aim now at providing an explicit form for (\ref{form}) by looking
for a suitable formula for the function $(\tau, t)\to
b_\varepsilon(\tau, t).$ To do this we substract (\ref{np}) where
$x(\tau)$ is replaced by $x_0(\tau+t)$ from (\ref{ps}) where
$x(\tau)$ is replaced by
$x_\varepsilon(\tau+t-\theta_\varepsilon(t))$ to obtain
\begin{eqnarray}\label{bis1}
& &\dot x_\varepsilon(\tau+t-\theta_\varepsilon(t))-\dot
  x_0(\tau+t)=\psi'(x_0(\tau+t))(x_\varepsilon(\tau+t-\theta_\varepsilon(t))-x_0(\tau+t))+\nonumber\\
&&+\>\varepsilon\phi(\tau+t-\theta_\varepsilon(t),x_\varepsilon(\tau+t-\theta_\varepsilon(t)),\varepsilon)
   + o(x_\varepsilon(\tau+t-\theta_\varepsilon(t))-x_0(\tau+t)).
\end{eqnarray}
By substituting (\ref{fw}) into (\ref{bis1}) and taking into account
that
\begin{eqnarray*}
& &\varepsilon a_\varepsilon(\tau,t)\psi'(x_0(\tau+t))\dot
x_0(\tau+t)=\varepsilon a_\varepsilon(\tau,t)\ddot x_0(\tau+t)\qquad
\mbox{and}\\
& &\varepsilon
b_\varepsilon(\tau,t)\psi'(x_0(\tau+t))y_1(\tau+t)=\varepsilon
b_\varepsilon(\tau,t)\dot y_1(\tau+t)
\end{eqnarray*}
we have
\begin{eqnarray*}
  \varepsilon\dot x_0(\tau+t)(a_\varepsilon)'_\tau(\tau,t)+\varepsilon y_1(\tau+t)
  (b_\varepsilon)'_\tau(\tau,t)&=&\varepsilon \phi(\tau+t-\theta_\varepsilon(t),x_\varepsilon(\tau+t-\theta_\varepsilon(t)))+\\
  &&+\> o(x_\varepsilon(\tau+t-\theta_\varepsilon(t))-x_0(\tau+t)),
\end{eqnarray*}
and so
\begin{eqnarray}
  \varepsilon
  (b_\varepsilon)'_\tau(\tau,t)&=&\varepsilon\left<z_1(\tau+t),\phi(\tau+t-\theta_\varepsilon(t)
  ,x_\varepsilon(\tau+t-\theta_\varepsilon(t)))\right>+\nonumber\\
  & &+\>\left<z_1(\tau+t),o(x_\varepsilon(\tau+t-\theta_\varepsilon(t))-x_0(\tau+t))\right>.\label{prob}
\end{eqnarray}
Moreover, since $z_1(\tau)=\rho_* z_1(\tau-T),$ from (\ref{b}) it
follows that
\begin{equation}\label{bc}
  b_\varepsilon(\tau,t)=\rho_* b_\varepsilon(\tau-T,t).
\end{equation}
System (\ref{prob})-(\ref{bc}) has a unique solution which, as it is
easy to verify, is given by the formula
\begin{eqnarray*}
  b_\varepsilon(\tau,t)& =&\frac{\rho_*}{\rho_*-1}\int_{\tau-T}^\tau \left<z_1(s+t),\phi(s+t-\theta_\varepsilon(t),
  x_\varepsilon(s+t-\theta_\varepsilon(t)),\varepsilon)\right>ds+\\
  &&+\>\frac{\rho_*}{\rho_*-1}\int_{\tau-T}^\tau\left<z_1(s+t),
  \frac{o(x_\varepsilon(s+t-\theta_\varepsilon(t))-x_0(s+t))}{\varepsilon}\right>ds.
\end{eqnarray*}
By substituting this formula into (\ref{form}) we obtain
\begin{eqnarray}
 & &\left< \dot
x_0(\tau+t)^\bot,x_\varepsilon(\tau+t-\theta_\varepsilon(t))-x_0(\tau+t)\right>
 = \label{perpe}\\
 &=&  \varepsilon k\frac{\|\dot
x_0(\tau+t)^\bot\|\rho_*}{\|z_1(\tau+t)\|(\rho_*-1)}\int_{\tau-T}^\tau
\left<z_1(s+t),\phi(s+t-\theta_\varepsilon(t),
  x_\varepsilon(s+t-\theta_\varepsilon(t)),\varepsilon)\right>ds+\nonumber \\
 & &+\> \varepsilon k\frac{\|\dot
x_0(\tau+t)^\bot\|\rho_*}{\|z_1(\tau+t)\|(\rho_*-1)}\int_{\tau-T}^\tau\left<z_1(s+t),\frac{o(x_\varepsilon(s+t-
\theta_\varepsilon(t))-x_0(s+t))}{\varepsilon}\right>ds,\nonumber
\end{eqnarray}
where $\dfrac{
o(x_\varepsilon(s+t-\theta_\varepsilon(t))-x_0(s+t))}{\varepsilon}\to
0$ as $\varepsilon\to 0$ uniformly in $s\in[-T,T]$ in virtue of
(\ref{ST}). On the other hand
\begin{eqnarray*}
& &\left< \dot
x_0(\tau+t)^\bot,x_\varepsilon(\tau+t-\theta_\varepsilon(t))-x_0(\tau+t)\right>=\\
& = & \left\|\dot
x_0(\tau+t)^\bot\right\|\left\|x_\varepsilon(\tau+t-\theta_\varepsilon(t))-x_0(\tau+t)
\right\|\,\cos\angle\left(\dot
x_0(\tau+t)^\bot,x_\varepsilon(\tau+t-\theta_\varepsilon(t))-x_0(\tau+t)\right)
\end{eqnarray*}
and by taking into account (\ref{st}) of Lemma 2,
(\ref{perpe}) and the fact that $\|\dot x_0(t)^\bot \|\not= 0$ for
any $t\in [0,T]$ we get
\begin{eqnarray*}
 \left\|x_\varepsilon(t-\theta_\varepsilon(t))-x_0(t)
\right\|&=&\varepsilon g(t) \int_{-T}^0
\left<z_1(s+t),\phi(s+t-\theta_\varepsilon(t),
  x_\varepsilon(s+t-\theta_\varepsilon(t)),\varepsilon)\right>ds+\\
  & &+\>\varepsilon g(t) \int_{-T}^0
\left<z_1(s+t),\dfrac{
o(x_\varepsilon(s+t-\theta_\varepsilon(t))-x_0(s+t))}{\varepsilon}\right>ds,
\end{eqnarray*}
where
$$
g(t)= \dfrac{k \rho_*}{\|z_1(t)\|(\rho_*-1)\,\cos\angle\left(\dot
x_0(t)^\bot,x_\varepsilon(t-\theta_\varepsilon(t))-x_0(t)\right)}
$$

\vskip0.2truecm\noindent is a continuous function on $[0,T]$ with
$g(t)\not= 0$ for any $t\in [0,T]$, Therefore,
\begin{eqnarray}\label{st11}
\left\|x_\varepsilon(t-\theta_\varepsilon(t))-x_0(t)
\right\|=\varepsilon g(t) \int_{-T}^0
\left<z_1(s+t),\phi(s+t-\theta_\varepsilon(t),
  x_\varepsilon(s+t-\theta_\varepsilon(t)),\varepsilon)\right>ds+o(\varepsilon).
\end{eqnarray}
On the other hand from Lemma 3 we have that $\Delta_\varepsilon(t)
\to 0$ uniformly in $t\in [0,T]$, thus we can rewrite (\ref{st11})
as follows
\begin{eqnarray*}
\left\|x_\varepsilon(t-\theta_\varepsilon(t))-x_0(t)
\right\|=\varepsilon g(t) \int_{-T}^0
\left<z_1(s+t),\phi(s+t-\theta_0,
  x_0(s+t),0)\right>ds+o(\varepsilon),
\end{eqnarray*}
introducing the change of variable $s+t=u$ in the integral we
finally get
\begin{equation*}
\left\|x_\varepsilon(t-\theta_\varepsilon(t))-x_0(t)
\right\|=\varepsilon g(t) \int_{t-T}^t \left<z_1(u),\phi(u-\theta_0,
  x_\varepsilon(u-\theta_0),0)\right>du+o(\varepsilon)
\end{equation*}
from which (\ref{conv2}) can be directly derived recalling that,
by Lemma~3, $x_\varepsilon(t-\theta_\varepsilon(t))\in
I(t,[-r_0,r_0])$ for any $\varepsilon\in (0,\varepsilon_0)$ and
any $t\in [0,T]$.

                                                        \qed

\vskip0.2cm\noindent As a straightforward consequence of Theorem 1
we have the following result.

\vskip0.2cm\noindent {\bf Corollary 1.} {\it Assume all the
conditions of Theorem~1, then for every $t\in[0,T]$ such that
$$
f_1(\theta_0,t)=0
$$
we have
$$
  \left\|x_\varepsilon(t-\theta_\varepsilon(t))-x_0(t)\right\|=o(\varepsilon),
$$
where $\theta_\varepsilon(t)\to \theta_0$ as $\varepsilon\to 0$ and
$x_\varepsilon(t-\theta_\varepsilon(t))\in I(t,[-r_0,r_0]).$}

\vskip0.3cm\noindent Next result is also a consequence of
Theorem~1.

\vskip0.2cm\noindent {\bf Corollary 2.} {\it Assume all the
conditions of Theorem~1. Moreover, assume that
\begin{equation}\label{pr11bis}
f_1(\theta_0,t)\not=0{\ \  for\ any\ \ } t\in[0,T].
\end{equation}
Then there exists $\varepsilon_1>0$ such that
$x_\varepsilon(s)\not=x_0(t)$ for any $s,t\in[0,T]$ and any
$\varepsilon\in(0,\varepsilon_1).$}

\noindent {\bf Proof.} Let $\varepsilon_0>0$ given by Theorem~1.
From (\ref{pr11bis}) we can choose
$\varepsilon_1\in(0,\varepsilon_0)$ in such a way that, for any
$\varepsilon\in(0,\varepsilon_1),$ we have both
\begin{equation}\label{forc}
o(\varepsilon)<\varepsilon M_1 |f_1(\theta_0,t)|,\quad{\rm for\ any\
}t\in[0,T],
\end{equation}
and the validity of (\ref{conv2}). Moreover, $\varepsilon_1$ can be
also chosen in such a way that there exists $\delta_0>0$ such that
the curve $\tau\to x_\varepsilon(\tau)$ intersects
$I(t,[-\delta_0,\delta_0])$ at only one point for any
$\varepsilon\in(0,\varepsilon_1)$ and $t\in[0,T].$ Such a choice is
possible, in fact, since $\dot x_{\varepsilon}(\tau-\theta_0)\to
\dot x_0(\tau)$ as $\varepsilon\to 0$ uniformly with respect to
$\tau\in[0,T]$ and the curve $r\to I(t,r)$ intersects the limit
cycle $x_0$ transversally at $r=0$ for any $t\in[0,T],$ then there
exists $\delta_0>0$ such that
$x_\varepsilon([t-\theta_0-\delta_0,t-\theta_0+\delta_0])$ and
$I(t,[-\delta_0,\delta_0])$ have only one common point for any
$t\in[0,T]$ and sufficiently small $\varepsilon>0.$ On the other
hand $\tau\to x_\varepsilon(\tau)$ cannot intersect
$I(t,[-\delta_0,\delta_0])$ for
$\tau\in\left[t-\theta_0-\frac{T}{2},t-\theta_0-\delta_0\right]
\cup\left[t-\theta_0+\delta_0,t-\theta_0+\frac{T}{2}\right]$ and
$\varepsilon>0$ sufficiently small, otherwise there would exist
sequences $\{\varepsilon_n\}_{n\in\mathbb{N}},$ $\varepsilon_n\to 0$
as $n\to\infty,$ $\{t_n\}_{n\in\mathbb{N}},$ $t_n\to t_0\in[0,T]$ as
$n\to\infty,$ $\{\tau_n\}_{n\in\mathbb{N}},$
$\tau_n\in[t_n-\theta_0-\frac{T}{2},
t_n-\theta_0-\delta_0]\cup[t_n-\theta_0+\delta_0,
t_n-\theta_0+\frac{T}{2}]$, $\tau_n\to \tau_0\in
\left[t_0-\theta_0-\frac{T}{2},t_0-\theta_0-\delta_0\right]\cup
\left[t_0-\theta_0+\delta_0,t_0-\theta_0+\frac{T}{2}\right]$ as
$n\to\infty$ such that  $x_{\varepsilon_n}(\tau_n)\in
I(t_n,[-\delta_0,\delta_0]),$ thus $x_0(\tau_0+\theta_0)=x_0(t_0),$
with $\tau_0+\theta_0\not= t_0$ and $|\tau_0+\theta_0-t_0|<T$, which
contradicts the fact that $T>0$ is the smallest period of $x_0.$

\noindent To conclude the proof assume now, by contradiction, that
there exist $\tilde\varepsilon\in(0,\varepsilon_1)$ and $\tilde s,
\tilde t\in [0,T]$ such that $x_{\tilde\varepsilon}(\tilde
s)=x_0(\tilde t).$ Since $\tau\to x_{\tilde\varepsilon}(\tau)$
intersects $I(\tilde t,[-\delta_0,\delta_0])$ at only one point
then Theorem~1 implies that $\tilde s=\tilde
t-\theta_{\tilde\varepsilon}(\tilde t).$ In conclusion, from
(\ref{conv2}) we have $\tilde\varepsilon M_1|f_1(\theta_0,\tilde
t)|\le o(\tilde\varepsilon)$ contradicting (\ref{forc}).

                                                               \qed

\vskip0.2cm\noindent The following result is crucial for the proof
of our existence result Theorem~3, but it can be also considered
as an independent contribution to the coincidence degree theory.

\vskip0.2cm\noindent {\bf Theorem 2.} {\it Assume conditions
(\ref{reg1}). For $\varepsilon>0$, let $G_\varepsilon :
C([0,T],\mathbb{R}^2)\to C([0,T],\mathbb{R}^2)$ be the operator
defined by
$$
(G_\varepsilon x)(t)=x(T)+\int_0^t
  \psi(x(\tau))d\tau+\varepsilon
  \int_0^t\phi(\tau,x(\tau),\varepsilon))d\tau, \qquad t\in[0,T].
$$
Let $W_{U_0}=\left\{x\in C([0,T],\mathbb{R}^2):x(t)\in
U_0,\right.$ for any $\left. t\in[0,T]\right\}.$ Assume that
\begin{equation}\label{sign}
\left<\dot x_0(0),z_0(0)\right>=\left< y_1(0),z_1(0)\right>=1.
\end{equation}
Finally, assume that for every $\theta_0\in[0,T]$ such that
$f_0(\theta_0)=0$ we have
\begin{equation}\label{pr11}
f_1(\theta_0,s+\theta_0)\not=0,\qquad{for\ any\ }s\in[0,T].
\end{equation}
Then, for all $\varepsilon>0$ sufficiently small,
$I-G_\varepsilon:C([0,T],\mathbb{R}^2)\to C([0,T],\mathbb{R}^2)$ is
not degenerate on the boundary of $W_{U_0}.$ Furthermore, there
exists a continuous vector field $F: \mathbb{R}^2\to \mathbb{R}^2$
such that
$$d(I-G_\varepsilon,W_{U_0})=d_B(F,U_0),$$
where $F(x_0(\theta))=f_0(\theta)\dot
x_0(\theta)+f_1(\theta,\theta)y_1(\theta)$ for any
$\theta\in[0,T].$}

\vskip0.2cm\noindent Some remarks are in order.

\vskip0.1cm\noindent {\bf Remark 1.} {\it As already observed
condition (\ref{sign}) does not affect the generality of
Theorem~2.}

\vskip0.2cm\noindent {\bf Remark 2.} {\it In Theorem~2 we could
replace $d_B(F,U_0)$ by $\,k\cdot {\rm ind}(x_0,F),$ where $k=+1$ or
$k=-1$ according with the orientation of the limit cycle $x_0.$
Precisely, $k=+1$ if the set $U_0$ is on the left side when one
follows $\partial U_0$ according to the parameterization $x_0(t)$
with $t$ increasing from $0$ to $T,$ and $k=-1$ in the opposite
case. Moreover, ${\rm ind}(x_0,F)$ is the Poincar\'e index of the
trajectory $x_0$ with respect to the vector field $F,$ namely the
total variation of an angle function of the vector $F(x_0(t))$ when
$t$ increases from $0$ to $T,$ see Lefschetz (\cite{lef}, Ch.~IX,
\S~4).}

\vskip0.2cm\noindent {\bf Remark 3.} {\it The Jordan theorem, see
Lefschetz (\cite{lef}, Theorem~4.7), ensures that the interior $U_0$
of $x_0$ does exist and it is an open set.}

\vskip0.3cm\noindent To prove Theorem 2 we need the following
preliminary lemma.

\vskip0.2cm\noindent {\bf Lemma 4.} {\it For any $s\in [0,T]$, let
\begin{equation}\label{Fs}
 F_s(\xi)=\int_{s-T}^s
 \Omega'_\xi(0,\tau,\Omega(\tau,0,\xi))\phi(\tau,\Omega(\tau,0,\xi),0)d\tau,
 \quad \mbox{for any} \;\;\; \xi\in \mathbb{R}^2.
\end{equation}
Then
\begin{eqnarray}
 & & \left<z_0(\theta),F_s(x_0(\theta))\right>=f_0(\theta){\ \  for\
  any}\ \ s,\theta\in[0,T],\nonumber\\
& & \left<z_1(\theta),F_s(x_0(\theta))\right>=f_1(\theta,s+\theta){\
\ for\
  any}\ \ s,\theta\in[0,T].\label{comp}
\end{eqnarray}
In particular,
\begin{equation}\label{Ftheta}
F_s(x_0(\theta))=\frac{1}{\left<\dot
x_0(t),z_0(t)\right>}f_0(\theta)\dot x_0(\theta)+\frac{1}{\left<
y_1(t),z_1(t)\right>}f_1(\theta,s+\theta)y_1(\theta)\quad { for\
any\ }s,\theta,t\in[0,T]. \end{equation} }

\noindent {\bf Proof.} It can be shown, see Krasnosel'skii
(\cite{kra}, Theorem~2.1), that
$\Omega'_\xi(t,0,x_0(\theta))=Y(t,\theta),$ where $Y(t,\theta)$ is
the fundamental matrix of the system
\begin{equation}\label{tp}
\dot y(t) = \psi'(x_0(t+\theta))y(t)
\end{equation}
satisfying $Y(0,\theta)=I$ and since
$\Omega'_\xi(0,t,\Omega(t,0,x_0(\theta)))\cdot\Omega'_\xi(t,0,x_0(\theta))=I$
we have
\begin{equation}\label{good}
  F_s(x_0(\theta))=
  \int_{s-T}^s
  Y^{-1}(\tau,\theta)\phi(\tau,x_0(\tau+\theta),0)d\tau.
\end{equation}
 Let us now show that
\begin{equation}\label{cl1}
  Y^{-1}(t,\theta)=Y(\theta,0) Y^{-1}(t+\theta,0).
\end{equation}
In fact, it is easy to see that $Y(t+\theta,0)$ is a fundamental
matrix for system (\ref{tp}) and so
$Y(t+\theta,0)Y^{-1}(\theta,0)$ is also a fundamental matrix for
(\ref{tp}), moreover we have that
$Y(t+\theta,0)Y^{-1}(\theta,0)=I$ at $t=0.$ Therefore
$Y(t+\theta,0)Y^{-1}(\theta,0)=Y(t,\theta)$ which is equivalent to
(\ref{cl1}).

\noindent By substituting (\ref{cl1}) into (\ref{good}) and by the
change of variable $\tau+\theta=t$ in the integral of (\ref{good})
we obtain
\begin{eqnarray*}
F_s(x_0(\theta))=
  Y(\theta,0)\int_{s-T}^{s}
  Y^{-1}(\tau+\theta,0)\phi(\tau,x_0(\tau+\theta),0)d\tau= Y(\theta,0)\int_{s-T+\theta}^{s+\theta}
  Y^{-1}(t,0)\phi(t-\theta,x_0(t),0)dt.
\end{eqnarray*}
\noindent Let $Z(t)$ be the fundamental matrix of system
(\ref{sopr}) given by $Z(t)=Z_0(t)Z_0^{-1}(0),$ where
$\,Z_0(t)=(z_0(t)\ z_1(t)), \newline t\in [0,T]$. Since
$Y^{-1}(t,0)=Z^*(t),$ see Perron \cite{perron} and Demidovich
(\cite{dem}, Sec. III, \S 12), then we have
$$
 F_s(x_0(\theta))=Y(\theta,0)\int_{s-T+\theta}^{s+\theta}
  Y^{-1}(\tau,0)\phi(\tau-\theta,x_0(\tau),0)d\tau
 =
\left(Z_0^*(\theta)\right)^{-1}\int_{s-T+\theta}^{s+\theta}
  Z_0^*(\tau)\phi(\tau-\theta,x_0(\tau),0)d\tau.
$$
Let
$$
 \Delta(s,\theta)=\int_{s+\theta-T}^{s+\theta} Z_0^*(\tau)
 \phi(\tau-\theta,x_0(\tau),0)d\tau,
$$
we have
\begin{eqnarray*}
& &\left<z_i(\theta),F_s(x_0(\theta))
\right>= \left<\left(\begin{array}{l} [z_i(\theta)]_1\\
{[z_i(\theta)]}_2\end{array}\right), {\left(\begin{array}{ll}
[z_0(\theta)]_1 & [z_0(\theta)]_2 \\
{[z_1(\theta)]}_1 & {[z_1(\theta)]}_2
\end{array}\right)}^{-1}\Delta(s,\theta)
 \right>=\nonumber\\
&=& \left<\left(\begin{array}{l} [z_i(\theta)]_1 \\
{[z_i(\theta)]}_2\end{array}\right),\ \frac{1}{{\rm
det}Z_0(\theta)} \left(\begin{array}{rr}
[z_1(\theta)]_2 & -[z_0(\theta)]_2 \\
{-[z_1(\theta)]}_1 & {[z_0(\theta)]}_1
\end{array}\right)\Delta(s,\theta)
 \right>=\nonumber\\
&=& \frac{1}{{\rm det}Z_0(\theta)}\left<\left(\begin{array}{l} [z_i(\theta)]_1 \\
{[z_i(\theta)]}_2\end{array}\right), \left(\begin{array}{ll}
[z_1(\theta)]_2 [\Delta(s,\theta)]_1 -[z_0(\theta)]_2 [\Delta(s,\theta)]_2\\
{-[z_1(\theta)]}_1 [\Delta(s,\theta)]_1 + {[z_0(\theta)]}_1
[\Delta(s,\theta)]_2
\end{array}\right)
 \right>=\nonumber\\
&=&   \frac{1}{{\rm det}Z_0(\theta)}\left\{ [z_i(\theta)]_1
[z_1(\theta)]_2 [\Delta(s,\theta)]_1 -[z_i(\theta)]_1
[z_0(\theta)]_2 [\Delta(s,\theta)]_2\right.-\nonumber\\
  & &\,- \left.[z_i(\theta)]_2 {[z_1(\theta)]}_1 [\Delta(s,\theta)]_1 +
[z_i(\theta)]_2 {[z_0(\theta)]}_1
[\Delta(s,\theta)]_2\right\}=\nonumber\\
&=&  \frac{1}{{\rm det}Z_0(\theta)}\left\{ [z_0(\theta)]_1
[z_1(\theta)]_2 -[z_0(\theta)]_2 {[z_1(\theta)]}_1
\right\}[\Delta(s,\theta)]_{i+1} =\nonumber\\
&=& \frac{1}{{\rm det}Z_0(\theta)}{\rm det}Z_0(\theta)
 \int_{s+\theta-T}^{s+\theta}
 \left<z_i(\tau),\phi(\tau-\theta,x_0(\tau),0)\right>d\tau.\nonumber
\end{eqnarray*}
For $i=0,1,$ we obtain (\ref{comp}). Furthermore, from
Lemma 1 and (\ref{sign}) we have that
$$
F_s(x_0(\theta))=\frac{1}{\left<\dot
x_0(t),z_0(t)\right>}f_0(\theta)\dot
x_0(\theta)+\frac{1}{\left<y_1(t),z_1(t)\right>}f_1(\theta,s+\theta)y_1(\theta)\quad
{\rm for\ any\ }\theta\in[0,T] \;\; {\rm and\ any\ } t\in[0,T].
$$

                                                                    \qed

\noindent {\bf Proof of Theorem 2.}

\noindent Let $\eta(t,s,\xi)$ be the solution of the system
\begin{equation}\label{aux}
  \dot q(t)=\psi'(\Omega(t,0,\xi))q(t)+\phi(t,\Omega(t,0,\xi))
\end{equation}
satisfying $\eta(s,s,\xi)=0$ whenever $\xi\in \mathbb{R}^2.$ It can
be shown, see (\cite{can}, Lemma~2), that
\begin{equation}\label{Feta}
F_s(\xi)=\eta(T,s,\xi)-\eta(0,s,\xi).
\end{equation}
Therefore, from (\ref{sign}), (\ref{pr11}) and (\ref{Ftheta}) we
have that
\begin{equation}\label{condeta}
  \eta(T,s,\xi)-\eta(0,s,\xi)\not=0
  {\ \ \rm\ \   for\ any\ \  } \xi\in\partial U_0{\ \ \rm\ \  and\ any\ \  }
  s\in[0,T],
\end{equation}
and by applying (\cite{non}, Theorem~2) we obtain the existence of
an $\varepsilon_0>0$ such that
$$
  d(I-G_\varepsilon, \widetilde{W}_{U_0})=d_B(\eta(T,0,\cdot), U_0) {\ \ \rm for\
  any}\ \ \varepsilon\in(0,\varepsilon_0),
$$
where
$
  \widetilde{W}_{U_0}=\{x\in {\rm C}([0,T],\mathbb{R}^2):\ \Omega(0,t,x(t))\in U_0\
  {\rm for\ any\ }
  t\in[0,T]\}.
$ Since as it is easy to see $\widetilde{W}_{U_0}=W_{U_0},$ by
taking into account (\ref{Ftheta}) and (\ref{Feta}) we end the proof
by defining $F(\xi)=F_0(\xi),\, \xi\in \mathbb{R}^2.$

                                                                   \qed

\vskip2mm\noindent The following existence theorem is the main
result of the paper.

\vskip0.2cm\noindent {\bf Theorem 3.} {\it Assume (\ref{reg1}).
Assume that for every zero $\theta_0\in[0,T]$ of the bifurcation
function $f_0$ we have
\begin{equation}\label{pr1}
f_1(\theta_0,t)\not=0{\ \  for\ any\ \ } t\in[0,T].
\end{equation}
Let $F\in C(\mathbb{R}^2,\mathbb{R}^2)$ be a vector field such that
on the boundary of $U_0$ it has the form
$F(x_0(\theta))=f_0(\theta)\dot
x_0(\theta)+f_1(\theta,\theta)y_1(\theta)$ for any $\theta\in[0,T].$
Assume
\begin{equation}\label{ad}
 d_B(F,U_0)\not=1.
\end{equation}
Then there exists $\varepsilon_0>0$ such that for every
$\varepsilon\in(0,\varepsilon_0)$ system (\ref{ps}) has at least
two $T$-periodic solutions $x_{1,\varepsilon}$ and
$x_{2,\varepsilon}$ satisfying
\begin{equation}\label{sat}
  x_{i,\varepsilon}(t-\theta_i)\to x_0(t)\quad {\rm as}\ \varepsilon\to
  0,\ i=1,2,
\end{equation}
where $\theta_1,\theta_2\in[0,T].$
 Moreover, we have that
$x_{1,\varepsilon}(t)\in U_0$ and $x_{2,\varepsilon}(t)\not\in U_0,$
for any $t\in[0,T]$ and any $\varepsilon\in(0,\varepsilon_0).$ }

\vskip3mm\noindent{\bf Proof.} \noindent Denote by
$W_\delta(\partial U_0)$ the $\delta$-neighborhood of the boundary
$\partial U_0$ of the set $U_0.$ Let $U_{1,\delta}=U_0\backslash
W_\delta(\partial U_0)$ and $U_{2,\delta}=U_0\cup
W_\delta(\partial U_0),$ thus the set $U_{1,\delta}$ tends to
$U_0$ from inside as $\delta\to 0,$ while $U_{2,\delta}$ tends to
$U_0$ from outside as $\delta\to 0.$ Since the limit cycle $x_0$
is isolated then there exists $\delta_0>0$ such that
\begin{equation}\label{defd}
G_0(x)\not=x {\quad \rm for\ any \quad  } x\in
\partial W_{U_{1,\delta}}\cup \partial W_{U_{2,\delta}}\quad \mbox{and any}\quad \delta\in (0, \delta_0].
\end{equation}
Moreover, being $T>0$, we can choose $\delta_0>0$ in such a way that
\begin{equation}\label{defdd}
\psi(\xi)\not=0 {\quad \rm for\ any \quad  } \xi \in\partial
U_{1,\delta}\cup \partial U_{2,\delta}\quad \rm and\ any\quad
\delta\in[0,\delta_0].
\end{equation}
From (\ref{defdd}) we get
$$
d_B(\psi,U_{1,\delta_0})=d_B(\psi,U_{2,\delta_0})=d_B(\psi,U_0).
$$
Since $U_0$ is the interior of the limit cycle $x_0$ of system
(\ref{np}) by Poincar\'e theorem, see Lefschetz (\cite{lef},
Theorem~11.1) or Krasnosel'skii et al. (\cite{per}, Theorem~2.3)
we have $d_B(\psi,U_0)=1$ and so
$$
d_B(\psi,U_{1,\delta_0})=d_B(\psi,U_{2,\delta_0})=1.
$$
In virtue of (\ref{defd}) and the fact that $W_U\cap\mathbb{R}^2=U,$
(\cite{maw}, Corollary~1) applies to conclude that
$$
d(I-G_0,W_{U_{1,\delta_0}})=d_B(\psi,U_{1,\delta_0}){\quad \rm and
\quad} d(I-G_0,W_{U_{2,\delta_0}})=d_B(\psi,U_{2,\delta_0}),
$$
hence
$$
d(I-G_0,W_{U_{1,\delta_0}})=1{\quad \rm and \quad}
d(I-G_0,W_{U_{2,\delta_0}})=1.
$$
Therefore, there exists $\varepsilon_0>0$ such that
\begin{equation}\label{f1}
d(I-G_\varepsilon,W_{U_{1,\delta_0}})=
d(I-G_\varepsilon,W_{U_{2,\delta_0}})=1\quad
  {\rm for\ any\ }\varepsilon\in(0,\varepsilon_0).
\end{equation}
Since by the definition of $z_1$ we have that
$z_1(t+T)=\rho_*z_1(t)$ for any $t\in[0,T],$ then for any
$t\in[0,T]$ it is easily seen that $f_1(\theta,t+T)=\rho_*
f_1(\theta,t)$, whenever $\theta\in[0,T],$ and thus from
(\ref{pr1}) we have also that $f_1(\theta_0,t+\theta_0)\not=0$ for
any $t\in[0,T].$ Therefore, all the conditions of Theorem~2 are
satisfied and we can take $\varepsilon_0>0$ sufficiently small to
have
\begin{equation}\label{f2}
  d(I-G_\varepsilon,W_{U_0})=d_B(F,U_0)\quad
  {\rm for\ any\ }\varepsilon\in(0,\varepsilon_0).
\end{equation}
By (\ref{ad}), (\ref{f1}) and (\ref{f2}) we conclude that for any
$\varepsilon\in(0,\varepsilon_0)$ there exist
\begin{equation}\label{sta2}
x_{1,\varepsilon}\in W_{U_0}\backslash W_{U_{1,\delta_0}},\ \ {\rm
and\ \ } x_{2,\varepsilon}\in W_{U_{2,\delta_0}} \backslash
W_{U_0}\end{equation} such that
$G_\varepsilon(x_{1,\varepsilon})=x_{1,\varepsilon}$ and
$G_\varepsilon(x_{2,\varepsilon})=x_{2,\varepsilon}.$ From
(\ref{sta2}) we have that for any $\varepsilon\in(0,\varepsilon_0)$
there exist points $t_{1,\varepsilon},t_{2,\varepsilon}\in[0,T]$
such that $x_{1,\varepsilon}(t_{1,\varepsilon})\in U_0\backslash
U_{1,\delta_0}$ and $x_{2,\varepsilon}(t_{2,\varepsilon})\in
U_{2,\delta_0}\backslash U_0.$ Thus $x_{1,\varepsilon}(t)\to\partial
U_0$ and $x_{2,\varepsilon}(t)\to\partial U_0$, for any $t\in
[0,T]$, as $\varepsilon\to 0,$ otherwise there would exist a
$T$-periodic solution $x_*$ to system (\ref{np}) and a point
$t_*\in[0,T]$ such that either $x_*(t_*)\in U_0\backslash
U_{1,\delta_0}$ or $x_*(t_*)\in U_{2,\delta_0}\backslash U_0$
contradicting (\ref{defd}). Therefore, see (\cite{malb},
Theorem~p.~287)  or (\cite{loud}, Lemma~2), for every
$i\in\left\{1,2\right\}$ there exists $\theta_i\in[0,T]$ satisfying
(\ref{sat}). The fact that $x_{1,\varepsilon}(t)\in U_0$ and
$x_{2,\varepsilon}(t)\not\in U_0$ for any $t\in [0,T]$ and
$\varepsilon>0$ sufficiently small follows from Corollary~2 and so
the proof is complete.

\qed

\noindent {\bf Remark 4.} {\it From the proof of Theorem~3 it
results that
$d(I-G_\varepsilon,W_{U_0}\backslash\overline{W}_{U_1,\delta_0})$
and $d(I-G_\varepsilon,\overline{W}_{U_2,\delta_0}\backslash
W_{U_0})$ are different from zero for
$\varepsilon\in(0,\varepsilon_0).$ This fact can be used to obtain
stability properties of solutions $x_{1,\varepsilon}$ and
$x_{2,\varepsilon}$ in the case when further information on the
number of $T$-periodic solutions to (\ref{ps}) belonging to  the
sets $W_{U_0}\backslash\overline{W}_{U_1,\delta_0}$ and
$\overline{W}_{U_2,\delta_0}\backslash W_{U_0})$ are available, see
Ortega \cite{ortega}.}

\

\noindent {\bf 3. An example.}

\noindent In this section we always assume that condition ($A_0$) is
satisfied. The well known formula by Poincar\'e, see Krasnoselskii
et. al. (\cite{per}, formula 1.16) states that
$$
 {\rm ind}(x_0,F)=\frac{1}{2\pi}\int_0^T\frac{[\alpha(\theta)]_1 [\alpha'(\theta)]_2-[\alpha(\theta)]_2
 [\alpha'(\theta)]_1}{[\alpha^2(\theta)]_1+[\alpha^2(\theta)]_2}d\theta,
$$
where $\alpha(\theta)=F(x_0(\theta)),$ $\theta\in[0,T].$ The
relationship between ${\rm ind}(x_0,F)$ and $d_B(F,U_0)$  was
discussed in Remark~2. In this section we show how the
representation $F(x_0(\theta))=f_0(\theta)\dot
x_0(\theta)+f_1(\theta,\theta)y_1(\theta), \theta\in[0,T],$ of the
function $F$ on $\partial U_0$ permits a simpler calculation of
$d_B(F,U_0).$ For this, we consider the case when $
\phi(t,\xi)=-\phi(t+T/2,\xi), $ which includes, in particular, the
classes of perturbations $\phi(t,x)=\sin t\cdot\phi_1(x)$ and
$\phi(t,x)=\cos t\cdot\phi_1(x),$ where $\phi_1\in
C(\mathbb{R}^2,\mathbb{R}^2).$ We can prove the following result.

\vskip0.1cm\noindent{\bf Proposition 1.} {\it Let $\widetilde{F}\in
C(\mathbb{R}^2,\mathbb{R}^2)$ be a vector field such that
$\widetilde{F}(x_0(\theta))=f_0(\theta)\dot
x_0(\theta)+f_1(\theta,T)\dot{x}_0^\bot(\theta),$ $\theta\in[0,T].$
Assume that
\begin{equation}\label{sign4}
\left<\dot x_0(0),z_0(0)\right>=\left< y_1(0),z_1(0)\right>=1,
\end{equation}
\begin{equation}\label{f0sym}
f_0(\theta)=-f_0(\theta+T/2) \quad { for\ any\ }\theta\in[0,T],
\end{equation}
\begin{equation}\label{f1sym}
f_1(\theta,T)=-f_1(\theta+T/2,T)\quad { for\ any\ }\theta\in[0,T].
\end{equation}
Moreover, assume that there exists an unique $\theta_0\in[0,T/2)$
such that $f_0(\theta_0)=0.$ Finally, assume that the function $f_0$
is strictly monotone at the point $\theta_0$ and that
\begin{equation}\label{pr111}
f_1(\theta_0,T)\not=0.
\end{equation}
Then either $d_B(\widetilde{F},U_0)=0$ or
$d_B(\widetilde{F},U_0)=2.$}

\vskip0.2cm

\noindent The proof of the proposition is based on the following
technical lemma.

\vskip0.1cm\noindent {\bf Lemma 5.} {\it Let $U\subset\mathbb{R}^2$
be an open set whose boundary $\partial U$ is a Jordan curve
$q:[0,T]\to\mathbb{R}^2,$ with $q(0)=q(T).$ Let
$\widetilde{F}:\mathbb{R}^2\to\mathbb{R}^2$ be a continuous vector
field such that $\widetilde{F}(\xi)\not=0$ for every $\xi\in\partial
U.$ Assume that for a continuous function $z:[0,T]\to\mathbb{R}^2,$
$z(0)=z(T),$ the following conditions hold:

1) $\left<z(\theta),\dot q(\theta)\right>\not=0$ for every
$\theta\in[0,T],$

2) the function
$f(\theta)=\left<z(\theta),\widetilde{F}(q(\theta))\right>$ has
exactly two zeros $\theta_1,\theta_2\in [0,T),$

3) the function $f$ is strictly monotone at $\theta_1$ and
$\theta_2,$

4) ${\rm sign}\left<z(\theta_1)^\bot,
\widetilde{F}(q(\theta_1))\right>=-{\rm sign}\left<z(\theta_2)^\bot,
\widetilde{F}(q(\theta_2))\right>.$

Then either $d_B(\widetilde{F},U)=0$ or $d_B(\widetilde{F},U)=2.$ }

\noindent{\bf Proof.} Assume that the parametrization $q$ is
positive, namely the set $U$ is on the left side if one follows
$\partial U$ according to the orientation given by $q$ when $t$
increases from $0$ to $T$, otherwise we consider the opposite
parametrization $\tilde q(\theta)=q(-\theta).$ For any $t\in[0,T]$
we denote by $\Theta(t)$ the angle (in radians) between the
vectors $\dot q(0)$ and $\dot q(t)$ calculated in the
counter-clockwise direction. Clearly $\Theta(t)$ is a multi-valued
function of $t.$ Let $\Gamma_{\dot q}(t)$ be the single-valued
branch of $\Theta(t)$ such that $\Gamma_{\dot q}(0)=0$ and let
$Q:\partial U\rightarrow \mathbb{R}^2$ be the vector field defined
by $Q(q(t)):=\dot q(t),$ whenever $t\in[0,T],$ hence $\Gamma_{\dot
q} (t)=\Gamma_{Q\circ q}(t).$ Following (\cite{per}, \S 1.2) the
function $t\rightarrow\Gamma_{\dot q}(t)$ is called the angle
function of the vector field $Q$ on the curve $q$ . Analogously,
considering the angle between $\widetilde{F}(q(0))$ and
$\widetilde{F}(q(t)),$ we can define the angle function
$\Gamma_{\widetilde{F}\circ q}(t)$ of the vector field
$\widetilde{F}$ on the curve $q.$ By the definition of the
rotation number for planar vector fields on the boundary of
simply-connected sets, see (\cite{per}, \S~1.3, formula~1.11) we
have
\begin{equation}\label{wehave}
d_B(\widetilde{F},U)=\frac{1}{2\pi}[\Gamma_{\widetilde{F}\circ
q}(T)-\Gamma_{\widetilde{F}\circ q}(0)].
\end{equation}
Therefore, in order to prove the lemma we must calculate the right
hand side of (\ref{wehave}). For this, denote by
$\widehat{h_1,h_2}\in[0,2\pi)$ the angle between the vectors $h_1$
and $h_2$ in the counter-clockwise direction, that is
$\widehat{h_1,h_2}+\widehat{h_2,h_1}=2\pi.$ Observe that
\begin{equation}\label{qqqq}
  \Gamma_{\widetilde{F}\circ q}(\theta)-\Gamma_{\dot
  q}(\theta)=\Gamma_{\dot q,\widetilde{F}\circ q}(\theta)-
      \widehat{\dot q(0),\widetilde{F}}(q(0)),
\end{equation}
where $\Gamma_{\dot q,\widetilde{F}\circ q}(\theta)$ is the single
valued branch of the multi-valued angle between $\dot q(\theta)$ and
$\widetilde{F}(q(\theta))$ such that $\Gamma_{\dot
q,\widetilde{F}\circ q}(0)=\widehat{\dot q(0),\widetilde{F}}(q(0)).$

To calculate $\Gamma_{\dot q,\widetilde{F}\circ q}(\theta)$ we
introduce the function
$\angle:\mathbb{R}^2\times\mathbb{R}^2\to[-\pi,\pi]$  as follows
$$
  \angle(h_1,h_2)=\left\{\begin{array}{ll} \widehat{h_1,h_2} & {\rm as\ \
  }\widehat{h_1,h_2}\in[0,\pi], \\  \widehat{h_1,h_2}-2\pi
  & {\rm as\ \
  }\widehat{h_1,h_2}\in(\pi,2\pi]\end{array}\right.
$$
By condition 3) we have that ${\rm ind}(\theta_i,f)=+1$ or ${\rm
ind}(\theta_i,f)=-1$ according to whether $f$ is increasing or
decreasing at $\theta_i,$ $i=1,2.$

Up to a shift in time, since $\theta_2-\theta_1<T$, we may assume
that the zeros $\theta_1,\theta_2$ of
$f(\theta)=\left<z(\theta),\widetilde{F}(q(\theta))\right>$ belong
to the interval $(0,T).$

Assume that $\left<z(\theta), \dot q(\theta)\right>>0$ for every
$\theta\in[0,T],$ otherwise we consider $\tilde
z(\theta)=z(-\theta)$ instead of $z(\theta).$ A possible way to
write explicity the function $\Gamma_{\dot q,\widetilde{F}\circ
q}(\theta)$ is the following
$$
  \Gamma_{\dot q,\widetilde{F}\circ q}(\theta)=
  \left\{
  \begin{array}{ll}
    \angle(z(\theta), \dot q(\theta))+\angle({\rm
   sign}\left<z(\theta),\widetilde{F}(q(\theta))\right>z(\theta),\widetilde{F}(q(\theta)))
   & \ \ {\rm as\ }\theta\in[0,\theta_1),\\
   \angle(z(\theta), \dot q(\theta))+\angle({\rm
   sign}\left<z(\theta),\widetilde{F}(q(\theta))\right>z(\theta),\widetilde{F}(q(\theta)))
   + & \\
   \hfill + \pi\,{\rm ind}(\theta_1,f){\rm sign}\left<z(\theta_1)^\bot,\widetilde{F}(q(\theta_1))\right>
    & \ \ {\rm as\ }\theta\in(\theta_1,\theta_2),\\
    \angle(z(\theta), \dot q(\theta))+\angle({\rm
   sign}\left<z(\theta),\widetilde{F}(q(\theta))\right>z(\theta),\widetilde{F}(q(\theta)))
   + & \\
   \hfill + \pi\,{\rm ind}(\theta_1,f){\rm
   sign}\left<z(\theta_1)^\bot,\widetilde{F}(q(\theta_1))\right>+&
   \\
   \hfill + \pi\,{\rm ind}(\theta_2,f){\rm sign}\left<z(\theta_2)^\bot,\widetilde{F}(q(\theta_2))\right>
    & \ \ {\rm as\ }\theta\in(\theta_2,T].
  \end{array}
  \right.
$$
It is easy to see that the above representation of the function
$\theta\rightarrow \Gamma_{\dot q,\widetilde{F}\circ q}(\theta)$
can be extend to $\theta_1$ and $\theta_2$ by continuity. Since
\begin{eqnarray*}
 \theta&\to&\angle(z(\theta),\dot q(\theta)),\\
 \theta&\to&\angle({\rm
   sign}\left<z(\theta),\widetilde{F}(q(\theta))\right>z(\theta),\widetilde{F}(q(\theta)))
\end{eqnarray*}
are $T$-periodic functions from (\ref{wehave})-(\ref{qqqq}), taking
into account that
\begin{equation}\label{qqq}
  d_B(Q,U)=\frac{1}{2\pi}[\Gamma_{\dot q}(T)-\Gamma_{\dot
  q}(0)]=1,
\end{equation}
(see e.g. \cite{per}, Theorem~2.4), we have
\begin{equation}
  d_B(\widetilde{F},
  U)= 1+\frac{1}{2}\left[{\rm ind}(\theta_1,f){\rm
sign}\left<z(\theta_1)^\bot,\widetilde{F}(q(\theta_1))\right>+
  {\rm ind}(\theta_2,f){\rm
  sign}\left<z(\theta_2)^\bot,\widetilde{F}(q(\theta_2))\right>\right].\label{lf}
\end{equation}
Since the function $f$ is $T$-periodic then
\begin{equation}\label{llf}
{\rm ind}(\theta_1,f)=-{\rm ind}(\theta_2,f)
\end{equation}
By assumption 4) and (\ref{llf}) the claim can be easily derived
from (\ref{lf}).

\qed

\noindent{\bf Proof of Proposition 1.}

\noindent Let $U=U_0,$  $q(t)=x_0(t),$
$z(t)=\dfrac{\dot{x}_0(t)}{\|\dot{x}_0(t)\|^2}, t\in [0,T],$ thus
the function $f_0$ turns out to be the function $f$ defined in
Lemma~5. Let us now show that all the conditions of Lemma~5 hold.
In fact, we have that $\left<\dot x_0(t),z(t)\right>= 1$ for any
$t\in [0,T]$ and so condition 1) is satisfied. Our assumptions
imply that the function $f_0$ has only two zeros
$\theta_1=\theta_0$ and $\theta_2=\theta_0+T/2$ in the interval
$[0,T]$ and it is strictly monotone at these points, thus
conditions 2) and 3) of Lemma~5 are also satisfied. Finally,
$\left<z(\theta)^\bot,
\widetilde{F}(x_0(\theta))\right>=f_1(\theta,T)$ and so
(\ref{f1sym})  implies condition 4) of Lemma~5. Hence the proof is
complete.

\qed

\noindent By combining Theorem 3 and Proposition 1 we obtain the
following result.

\vskip0.1cm\noindent{\bf Corollary 3.} {\it Assume conditions
(\ref{reg1}) and assume that
$$
\phi(t,\xi)=-\phi(t+T/2,\xi)\quad{ for\ any\ }t\in[0,T]\ {\it\ and \
any\ } \xi\in\mathbb{R}^2.
$$
 \noindent Moreover, assume that there exists a unique
$\theta_0\in[0,T/2)$ such that $f_0(\theta_0)=0.$ Finally, assume
that the function $f_0$ is strictly monotone at the point $\theta_0$
and
\begin{equation}\label{52}
f_1(\theta_0,t)\not=0\quad{ for\ any\ }t\in[0,T].
\end{equation}
Then there exists $\varepsilon_0>0$ such that for every
$\varepsilon\in(0,\varepsilon_0)$ system (\ref{ps}) has at least
two $T$-periodic solutions $x_{1,\varepsilon}$ and
$x_{2,\varepsilon}$ satisfying
$$
  x_{i,\varepsilon}(t-\theta_i)\to x_0(t)\quad { as}\ \varepsilon\to
  0,\ i=1,2,
$$
where  $\theta_1,\theta_2\in\left\{\theta_0,\theta_0+T/2\right\}.$
Furthermore, we have that $x_{1,\varepsilon}(t)\in U_0$ and
$x_{2,\varepsilon}(t)\not\in U_0,$ for every $t\in[0,T]$ and
$\varepsilon\in(0,\varepsilon_0).$ }

\vskip0.2cm\noindent {\bf Proof of Corollary 3.}

\noindent To apply Theorem~3 we only have to verify condition
(\ref{ad}). For this we will make use of Proposition~1. Without
loss of generality we can assume that
\begin{equation}\label{wo}
  \left<y_1(0),\dot{x}^\bot(0)\right>>0.
\end{equation}
We claim that, under the conditions of Corollary~3,  the vector
field $\widetilde{F}$ of Proposition~1 is homotopic on $\partial
U_0$ to the vector field $F$ of Theorem~3. To prove the claim we
show that the following homotopy joining $\widetilde{F}$ and $F$
$$
  D_\lambda(x_0(\theta))=f_0(\theta)\dot{x}_0(\theta)+f_1(\theta,\lambda
  T+(1-\lambda)\theta)(\lambda\dot{x}_0(\theta)^\bot+(1-\lambda)y_1(\theta)),\quad \mbox{with}\;
  \theta\in[0,T]\quad\mbox{and}\;
  \lambda\in[0,1],
$$
is admissible. Assume the contrary, therefore there exist
$\lambda_0\in[0,1]$ and $\theta_0\in[0,T]$ such that
$$
  f_0(\theta_0)\dot{x}_0(\theta_0)+f_1(\theta_0,\lambda_0
  T+(1-\lambda_0)\theta_0)(\lambda_0\dot{x}_0(\theta_0)^\bot+(1-\lambda_0)y_1(\theta_0))=0.
$$
By condition (\ref{wo}) we have that the vectors
$\dot{x}_0(\theta_0)$ and
$\lambda_0\dot{x}_0(\theta_0)^\bot+(1-\lambda_0)y_1(\theta_0)$ are
linearly independent thus
$$
  f_0(\theta_0)=0\quad\mbox{and}\quad f_1(\theta_0,\lambda_0
  T+(1-\lambda_0)\theta_0)=0
$$
contradicting assumption (\ref{52}). Hence we have proved that
$$
  d_B(F,U_0)=d_B(\widetilde{F},U_0).
$$
Applying Proposition~1 we obtain that
$$
  d_B(F,U_0)\in\{0,2\},
$$
namely assumption (\ref{ad}) of Theorem~3 is satisfied and the
conclusion of the corollary follows from Theorem~3.

\qed

\vskip0.4truecm\noindent At the end of the paper we would like to
stress that all the functions $y_1,$ $z_0$ and $z_1$ can be easily
determined both analytically and numerically once the limit cycle
$x_0$ is known. We give in the following a sketch of both
approaches.

\vskip0.4truecm\noindent 1) The analytical approach.

\noindent Since $\dot x_0$ is one of the two eigenfunctions of
system (\ref{ls}) then by using well known formulas, see e.g.
Pontrjagin (\cite{pon}, p.~138), the dimension of the system
(\ref{ls}) can be decreased by 1, thus the obtained
one-dimensional system can be easily solved to determine $y_1.$
Furthermore, by Lemma~1 the eigenfunctions $z_0$ and $z_1$ can be
determined by the formula
$$
  (z_0(t)\ z_1(t))=\left((\dot x_0(t)\ y_1(t))^*\right)^{-1}.
$$

\vskip0.4truecm\noindent 2) A direct numerical approach.

\noindent From Lemma 1 we have $\left<\dot x_0(0),z_1(0)\right>=0,$
therefore as initial condition we may take $z_1(0)=\dot x_0(0)^\bot$
and then $z_1$ can be obtained by a numerical computation. By the
definition of $z_1$ there exists a $T$-periodic function $a\in
C(\mathbb{R},\mathbb{R}^2)$ such that $z_1(t)=a(t)\,{\rm e}^{\rho_*
t}.$ Assume, that $\rho_*<0.$ Let us fix an arbitrary vector
$\xi\in\mathbb{R}^2,$ which is linearly dependent with $z_1(0)$ and
calculate the solution $z$ of system (\ref{sopr}) satisfying
$z(0)=\xi$ on the interval $[0,kT]$ where $k\in\mathbb{N}.$ It turns
out that larger is $k$ better accuracy is obtained.  Observe that
$z$ can be represent by
\begin{equation}\label{z}
  z(t)=\alpha a(t){\rm e}^{\rho_* t}+z_0(t),
\end{equation}
where $z_0$ is an eigenfunction of (\ref{sopr}), and since ${\rm
e}^{\rho_* t}\to 0$ as $t\to+\infty$ then for given $k\in\mathbb{N}$
we may take
$$
  z_0(t)=z(t+(k-1)T)\quad{\rm for\ any\ }t\in[0,T].
$$
For the case $\rho_*>0$ one should make the change of variables
$\tilde z(t)=z(-t)$ for any $t\in\mathbb{R}$ in (\ref{sopr}), to
calculate $\tilde z_0$ on $[-kT,0]$ (for this it is necessary to
expand, the function $z_1$ on the interval $[-kT,0]$) and then put
$z_0(t)=z_1(-t)$ for any $t\in[0,T].$ Once $z_0$ is calculated
with the desirable accuracy the function $y_1$ can be determined
as the solution of (\ref{ls}) with initial condition
$y_1(0)=z_0(0)^\bot.$

\end{document}